\documentclass[11pt,leqno]{amsart}

\usepackage[utf8x]{inputenc} 	        
\usepackage{lmodern}                  
\usepackage[T1]{fontenc}              
\usepackage[english]{babel}           
\usepackage{ucs}             	        
\usepackage{hyperref}

\usepackage{geometry}
\usepackage{graphicx}     
\usepackage{xcolor}       
\usepackage{stmaryrd}

\usepackage{amsmath, amsthm}
\usepackage{amsfonts, amssymb}
\usepackage{bm}                 

\theoremstyle{plain}
	\newtheorem{theorem}{Theorem}[section]
	\newtheorem*{theorem*}{Theorem}
	
	\newtheorem{proposition}{Proposition}[section]
	\newtheorem*{proposition*}{Proposition}

\numberwithin{equation}{section}

\theoremstyle{remark}
	\newtheorem{remark}{\textbf{Remark}}[]

\theoremstyle{definition}

\DeclareMathOperator{\supp}{supp}

\DeclareMathOperator{\mes}{mes}

\newcommand{\eg}{\emph{e.g. }}

\newcommand{\RR}{\mathbb{R}}

\newcommand{\ZZ}{\mathbb{Z}}
\newcommand{\TT}{\mathbb{T}}

\newcommand{\Q}{\mathcal{Q}}
\newcommand{\AAA}{\mathcal{A}}
\newcommand{\BBB}{\mathcal{B}}

\newcommand{\grad}{\nabla}
\newcommand{\ds}{\displaystyle}

\def\R{\mathbb{R}}

\newcommand{\upwind}{\textbf{upwind}}
\newcommand{\mcu}[1]{\textbf{MCU(#1)}}
\newcommand{\mcua}{\mcu{1}}
\newcommand{\mcub}{\mcu{0}}

\title{An Exact Rescaling Velocity Method for some Kinetic Flocking Models}

\author[T. Rey]{Thomas Rey}
\address{Thomas Rey \\
Laboratoire P. Painlevé, CNRS UMR 8524 \\
Université Lille 1 \\
59655 Villeneuve d'Ascq Cedex\\
France
}
\email{thomas.rey@math.univ-lille1.fr}

\author[C. Tan]{Changhui Tan}
\address{Changhui Tan \\
Center of Scientific Computation and Mathematical Modeling (CSCAMM) \\
The University of Maryland \\
College Park, MD, 20742-4015\\
USA}
\email{ctan@cscamm.umd.edu}

\keywords{flocking models, kinetic equations, rescaling velocity methods, finite volume methods, upwind scheme, large time behavior}
	
\subjclass[2010]{Primary:   82C40, 
                 Secondary: 65N08, 
}

\date{Final version}

\begin{document}

  \begin{abstract}
    In this work, we discuss kinetic descriptions of flocking models, of the so-called Cucker-Smale \cite{CuckerSmale:2007} and Motsch-Tadmor \cite{MotschTadmor:2011} types. 
    These models are given by Vlasov-type equations where the interactions taken into account are only given long-range bi-particles interaction potential.
    We introduce a new exact rescaling velocity method, inspired by the recent work \cite{FilbetRey:2012}, allowing to observe numerically the flocking behavior of the solutions to these equations, without a need of remeshing or taking a very fine grid in the velocity space.
    To stabilize the exact method, we also introduce a modification of the classical upwind finite volume scheme which preserves the physical properties of the solution, such as momentum conservation. 
  \end{abstract}
  
  \maketitle
  
  \tableofcontents

	\section{Kinetic Description of Flocking Models}
	  \label{secKinDescript}
	  	  
		We are interested in this paper with numerical simulations of Vlasov-type kinetic description of flocking models
		\begin{equation}\label{eq:flocking}
		  \left \{ \begin{aligned}
	  	  & \frac{\partial f}{\partial t}+v\cdot\grad_x f+\grad_v\cdot \Q(f)=0, \ \forall x \in \Omega \subset \RR^d, \, v \in \R^d, \\
	  	  &  \, \\
	  	  & f(0,x,v) = f_0(x,v).
	  	\end{aligned} \right.
		\end{equation}
		The \emph{distribution function} $f=f(t,x,v)$ describes the probability to find an individual at time $t >0$ at the infinitesimal position of the phase space $dx \, dv$.
		The set $\Omega$ will be either $\TT^d$ or $\RR^d$.
		The integral operator $\Q$, the \emph{flocking} operator, characterizes the nonlocal \emph{interactions}. Typical examples are given by the so-called \emph{Cucker-Smale} \cite{CuckerSmale:2007} model and \emph{Motsch-Tadmor}
		\cite{MotschTadmor:2011} model, where the operators (See \eg the paper from Tadmor and Ha \cite{TadmorHa:2008} for details about the derivation) are given by
		\begin{align}
	  	\label{eq:collisionCS} \text{Cucker-Smale model: } & \Q_{CS}(f)=\iint\phi(|x-y|)(v^*-v)f(y,v^*)f(x,v) \,dv^* \,dy,\\
		  \label{eq:collisionMT} \text{Motsch-Tadmor model: } & \Q_{MT}(f)=\left[\iint\phi(|x-y|)f(y,v^*) \, dv^* \, dy\right]^{-1}Q_{CS}(f).
		\end{align}

		The function $\phi$ is the \emph{influence function}, and characterizes the range of the interactions between individuals.
		If $\phi$ decays slowly enough at infinity, the
		system converges to a \emph{flock}, namely all the individuals travel in a close packed regime, at a constant speed. 
		In fact, if one defines
		\[
		  S(t):=\sup_{(x,v),(y,v^*)\in \supp(f(t))}|x-y|,\quad
		  V(t):=\sup_{(x,v),(y,v^*)\in \supp(f(t))}|v-v^*|
		\]
		to be respectively the largest variation in position and in velocity for the system, one can prove the
		following theorem:

		\begin{theorem}[\cite{Carrillo:2010,Tan:2014}]\label{thm:flocking}
  		Suppose $\ds\int^\infty\phi(r)dr=\infty$, then
	  	$S(t)$ is bounded for all times, and $V(t)$ decays to 0 exponentially in time.
		\end{theorem}

		From this theorem, we know that  at least for the two flocking operators of interest, the equilibrium states $f_\infty$ of the system are \emph{monokinetic}: they
		have the form 
		\[
		  f_\infty(x,v) = \rho(x) \delta(v-\bm{c}),
		\]
		where $\bm{c}\in \RR^d$ is a constant velocity which
                depends on the initial condition\footnote{If the
                  flocking operator is symmetric, as in the
                  Cucker-Smale case, this quantity is given by the
                  initial average velocity of the system.} and $\rho$ is the macroscopic \emph{density} of the system:
		\[
		  \rho(x) = \int_{\RR^d} f(x,v) \, dv.
		\]
		With such exponentially fast creation of $\delta$-singularities, when designing a numerical method for \eqref{eq:flocking}, one cannot expect to achieve a correct accuracy for large time (and this is particularly true when using a high accuracy spectral method\footnote{This is a natural choice of approximation of the flocking operator, because of its convolution structure.}, because of Gibbs phenomenon \cite{canuto:88}). It is then challenging to design a numerical scheme which captures the blow-up correctly.

		In \cite{Tan:2014}, schemes based on discontinuous Galerkin method are
		derived to deal with $\delta$-singularities (flocking and
		clustering). In this paper,  as we mentioned earlier, we shall only focus
		on the flocking case. 
		To this end, we shall introduce a technique based on the information provided by the hydrodynamic fields computed from a macroscopic model corresponding to the original kinetic equation. Then by rescaling the kinetic equation using the knowledge of its qualitative behavior (mostly Theorem \ref{thm:flocking}), we will solve another equation which does not exhibit concentration. The original solution will finally be obtained by reverting the rescaling.

    Recently, F. Filbet and G. Russo proposed in~\cite{FilbetRusso:2004} a rescaling method for space homogeneous kinetic equations on a fixed grid. This idea is mainly based on the self-similar behavior of the solution to the kinetic equation. However for spatially inhomogeneous case, the situation is much more complicated and this method have not been applied since the transport operator and the boundary conditions break down this self-similar behavior. Then, F. Filbet and the first author proposed in \cite{FilbetRey:2012} an extension of this method to the space inhomogeneous case using an approximate closure of the macroscopic equations based on the knowledge of the hydrodynamic limit of the system. This is particularly well suited to the study of granular media.

		In this work, we shall come back to the original, exact approach, and give a new method to couple the evolution of the kinetic equation together with the computation of the rescaling function. It is based on a new modification of the classical upwind fluxes (see \eg \cite{leveque:2002} for a complete introduction on the topic) allowing to take into account some of the physical properties of the equation.
		To follow the flock during time, we will also introduce a shift in velocity for the rescaling function.
		Besides capturing the correct flocking behavior of the model, the rescaled equation is also local in the velocity variable, and hence computationally cheaper than the original nonlocal equation for $f$.

  \section{Scaling on Velocity}
    \label{secScaling}
    
    This section is devoted to the presentation of the  scaling for equation~\eqref{eq:flocking} allowing to follow the change of scales in velocity. 
    It is an extension to the space-dependent setting of the method first introduced in \cite{FilbetRusso:2004}, using the relative kinetic energy as a scaling function. 
    This method was reminiscent from the work of A. Bobylev, J.A. Carrillo, and I. Gamba  \cite{bobylev:2000} about Enskog-like inelastic interactions models. Indeed, in one section of this work, the authors scaled the solution of the spatially homogeneous collision equation by its thermal velocity, in order to study a drift-collision equation, where no blow-up occurs.
    The same technique was also used by  S. Mischler and C. Mouhot in \cite{Mischler:2006} to prove the existence of self-similar solutions to the granular gases equation.

    For a given positive function $\omega:\RR^+\times\Omega\mapsto \RR^+$, we introduce a new distribution $ g (t, x, \xi)$ by setting
		\begin{equation} 
		  \label{def:Rescaling}
	  	f(t,x,v)=\omega(t,x)^d g(t,x,\xi),\quad \xi=\omega(v-\bm{u}),
		\end{equation}		
		where the function $\omega$ (or more precisely $\omega^{-1}$), the \emph{scaling factor}, is assumed to be an accurate measure of the ``support'' or  scale  of the distribution $f$ in velocity variables. Then according to this scaling, the distribution $g $ should naturally ``follow'' either the concentration or the spreading in velocity of the distribution $f$. 
		
		Moreover, it is straightforward using \eqref{def:Rescaling} to see that $g$
		has the following qualitative properties:
		\begin{enumerate}
		\item Its local density is the same as the one of the original distribution:
		  \begin{equation}
		    \label{eq:scalingDensity}
		    \rho(t,x):= \int_{\RR^d} g(t,x,\xi) \, d\xi=\int_{\RR^d} f(t,x,v) \, dv.
		  \end{equation}
		\item Due to the shift in velocity, its local momentum is everywhere $0$:
		  \begin{equation}
		    \label{eq:scalingMomentum}
		    M(t,x):= \int_{\RR^d} g(t,x,\xi) \, \xi \, d\xi=0, \ \forall t\geq 0, \, x\in \Omega.
		  \end{equation}
		\end{enumerate}

		The question now is to find an appropriate scaling factor $\omega$ so that $g$
		neither vanishes nor becomes singular in all time.
		
		\begin{remark}
	  	From Theorem \ref{thm:flocking}, $V(t)$ is decaying exponentially fast in
	  	time. Since this quantity is essentially the support of $f$, in order for the support of $g$ to remain bounded, we expect that $\omega$ should grow exponentially in time.
		\end{remark}
				
		\subsection{A Spatially ``Homogeneous'' System}
		  \label{subHomog}
  		
  		We first consider the dynamics without the free transport term:
	  	\begin{equation}\label{eq:flock}
		    \frac{\partial f}{\partial t}+\grad_v\cdot \Q(f)=0,
		  \end{equation}
                  since the flocking operator $\Q$ is the main driving
                  force towards velocity concentration. 
		  Note that here, the system is not completely
                  spatially homogeneous, as $\Q$ is a nonlocal
                  operator in space.

			Plugging the expression of the scaling function \eqref{def:Rescaling} into the flocking equation \eqref{eq:flock}, we have for the first term	
			\begin{align*}
				\partial_t f~=~&d\omega^{d-1}\partial_t\omega g+\omega^d\big[\partial_t g+\grad_\xi
				  g\cdot\left(\partial_t \omega(v-\bm{u})-\omega \partial_t\bm{u}\right)\big]\\
				 =~&\omega^d\left[\partial_t g+\frac{\partial_t\omega}{\omega}\grad_\xi\cdot(\xi
				g)\right]-\omega^{d+1}\partial_t\bm{u}\cdot\grad_\xi g.
		  \end{align*}
		  Moreover, concerning the flocking operator, one has to distinguish between Cucker-Smale and Motsch-Tadmor. The former one yields
		  \begin{align*}
				\Q_{CS}(f)~=~&\omega^d\iint\phi(|x-y|)\left(\frac{\xi^*}{\omega(y)}+\bm{u}(y)-\frac{\xi}{\omega(x)}-\bm{u}(x)\right)g(x,\xi)\, g(y,\xi^*) \,d\xi^*\,
				dy\\
				~=~&\omega^dg(x,\xi)\int\phi(|x-y|)(\bm{u}(y)-\bm{u}(x))\rho(y)\,dy-\omega^{d-1}
				\xi \,g(x,\xi)\int\phi(|x-y|)\rho(y)\,dy,
			\end{align*}
			whereas the latter one yields
			\begin{align*}
				\Q_{MT}(f)~=~&\left[\int\phi(|x-y|)\,\rho(y)\,dy\right]^{-1}\Q_{CS}(f)\\
				~=~&\omega^dg(x,\xi) \frac{\int\phi(|x-y|)(\bm{u}(y)-\bm{u}(x))\rho(y)\,dy}{\int\phi(|x-y|)\rho(y)\,dy}-\omega^{d-1}
				\xi \, g(x,\xi).
			\end{align*}
      Gathering everything and  using the chain rule, we obtain the following equation for $g$:
      \begin{equation}
        \label{eq:gGeneralOmega}
        \frac{\partial g}{\partial t}+\left[\frac{\partial_t\omega}{\omega}-
        \AAA(t,x)\right]\grad_\xi\cdot(\xi g)-\omega\left[\partial_t\bm{u}-\BBB(t,x)\right]\cdot\grad_\xi
        g=0,
      \end{equation}
      where the operators $\AAA$ and $\BBB$ are functions of the macroscopic quantities and the influence function only, and depend on the model considered. 
      More precisely, we have for the Cucker-Smale model \eqref{eq:collisionCS}
      \begin{equation}
        \label{eq:AB_CS}
        \AAA_{CS}(t,x) := \int\phi(|x-y|)\rho(t,y)\,dy,\quad \BBB_{CS}(t,x):=\int\phi(|x-y|)(\bm{u}(t,y)-\bm{u}(t,x))\rho(t,y)\,dy,
      \end{equation}
      and for the Motsch-Tadmor model \eqref{eq:collisionMT}
      \begin{equation}
        \label{eq:AB_MT}
        \AAA_{MT}(t,x)=1, \quad \BBB_{MT}(x)=\frac{\int\phi(|x-y|)(\bm{u}(t,y)-\bm{u}(t,x))\rho(t,y)\,dy}{\int\phi(|x-y|)\rho(t,y)\,dy}.
      \end{equation}

      Computing the zeroth and first moments in velocity of equation \eqref{eq:flock}, we get the following
evolutions for the macroscopic quantities: 
			\begin{align*}
  			\frac{d}{dt}\rho(t,x)&=\frac{d}{dt}\int f(t,x,v)\,dv=0,\\
	  		\frac{d}{dt}(\rho(t,x)\bm{u}(t,x))&=\frac{d}{dt}\int f(t,x,v)\,v\,dv=\int \Q(f) \, dv=\rho(t,x)\BBB(t,x).
			\end{align*}
			It implies that the mass $\rho(t,x)=\rho(x)$ is constant in time, and that the following evolution law holds
			\begin{equation}
			  \label{eq:evolUHomog}
			  \left \{ \begin{aligned}
   			  & \partial_t\bm{u} -\BBB(t,x)=0, \\
   			  & \, \\
		  	  & \bm{u}(0,x) = \frac{1}{\rho(x)} \int_{\RR^d} f_0(x,v) \,v \, dv.
			  \end{aligned} \right.  
			\end{equation}
			Moreover, $\AAA$ is independent in time for both models, since it does not depend on $\bm{u}$.

      We have now enough information to define the scaling function. Let us set
      \begin{equation}
        \label{def:omega}
        \omega(t,x)=\omega_0(x)\exp\left[t \AAA(x)\right], \ \forall t \geq 0, \ x \in \Omega,
      \end{equation}
      for any measurable, positive function $\omega_0$. Note that such an $\omega$ behaves as expected, namely grows exponentially in time to compensate the concentration, as $\AAA$ is nonnegative. This is particularly true for the Motsch-Tadmor model \eqref{eq:collisionMT}, where we have according to \eqref{eq:AB_MT} the explicit form
      \begin{equation}
        \label{eq:omegaMT} 
        \omega(t,x) = \omega_0(x) \, e^t, \ \forall t \geq 0, \ x \in \Omega.
      \end{equation}
      
      Then, plugging both \eqref{eq:evolUHomog} and \eqref{def:omega} in \eqref{eq:gGeneralOmega}, one obtains that
      \[\frac{\partial g}{\partial t}=0.\] 
      This provides a perfect scaling of the system, and since $g$ is
      constant in time, the whole dynamics of $f$ is given by the
      dynamics of the scaling function $\omega$ and the macroscopic
      velocity $\bm{u}$. 
      More precisely, we have for $f$ using \eqref{def:Rescaling}
			\begin{align*}
				f(t,x,v)=& ~\omega(t,x)^d g\left (t,x,\omega(t,x)(v-\bm{u}(t,x))\right )\\
				=&~\omega(t,x)^d g(0,x,\omega(t,x)(v-\bm{u}(t,x)))\\
				=&~\frac{\omega(t,x)^d}{\omega_0(x)^d}
				f\left(0,x,\frac{\omega(t,x)}{\omega_0(x)}(v-\bm{u}(t,x))+\bm{u}(0,x)\right)\\
				=&~\exp \left [d\,t\AAA(x) \right]f_0\left(x,\exp \left [t\AAA(x)\right ]v+\bm{u}(0,x)-\exp \left [t\AAA(x)\right ]\bm{u}(t,x)\right),
			\end{align*}
			the momentum $\bm u$ being a solution to \eqref{eq:evolUHomog}. Such an $f$ is a solution to \eqref{eq:flock}. 
			
			\begin{remark}
			  Note that the integral term $\BBB$ has a convolution structure in $x$, and a spectral method could be used to propagate $\bm{u}$ with no difficulty and high accuracy.
      \end{remark}

    \subsection{The Full Rescaled System}
      \label{subFull}
      
      We are now ready to go back to the full system \eqref{eq:flocking}.  Applying the scaling $\omega$ introduced in the last section, we have that:
      \[
				v\cdot\grad_xf~=~\left(\frac{\xi}{\omega}+\bm{u}\right)\cdot\left[
				d\omega^{d-1}\grad_x\omega~g+\omega^d\left(\grad_xg+
				\frac{\grad_x\omega}{\omega}~(\xi\cdot\grad_x)g-\omega\sum_i\partial_{\xi_i}g\grad_xu_i\right)\right].
		  \]
		  With the transport term, the rescaled $g$ will not be constant anymore as time goes by. 
		  A direct computation using \eqref{eq:gGeneralOmega} yields the following dynamics for $g$:
			\begin{align}
			  \label{eq:gNonConservative}
				\frac{\partial g}{\partial t} 
				 & +\left[\frac{\partial_t \omega}{\omega}+\bm{u}\cdot\frac{\grad_x\omega}{\omega}-\AAA(t,x)\right]\grad_\xi\cdot (\xi g)-\omega\left[\partial_t \bm{u} +\bm{u}\cdot\grad_x \bm{u}-\BBB(t,x)\right]\cdot\grad_\xi g + \bm{u}\cdot\grad_xg \\
				 & +\frac{\xi}{\omega}\cdot\left[\frac{\grad_x\omega}{\omega}\grad_\xi\cdot(\xi g) + \grad_xg - \omega\sum_i\partial_{\xi_i}g\grad_xu_i\right]=0. \notag
			\end{align}
			After some computations, one can rewrite this equation on the following form:
			\begin{align}
			  \label{eq:gAlmostConservative}
				\frac{\partial g}{\partial t} 
				 & + \nabla_\xi \cdot \left \{  \left [ -	\omega \left( \partial_t \bm{u} + \bm{u} \cdot \nabla_x \bm{u} - \BBB \right ) 
				   - \xi\cdot\grad_x\bm{u} 
				   + \frac{1}{\omega} \left( \partial_t \omega - \omega \AAA + \left(\bm{u}+\frac{\xi}{\omega}\right) \cdot \nabla_x \omega \right) \xi
				   \right ] g \right \} \\
				 & + \nabla_x \cdot \left [ \left (\bm{u} + \frac{\xi}{\omega}\right ) g\right]	= 0. \notag
		  \end{align}
		  
		  Multiplying the original flocking equation \eqref{eq:flocking} by respectively $1$ and $v$ and integrating in the velocity variable, we obtain by using the definition of $g$ \eqref{def:Rescaling} the evolution of the macroscopic quantities:
      \begin{equation}
        \label{eq:sysMomentCoupled}			
			  \left \{\begin{aligned}
			   & \partial_t \rho+\grad_x\cdot(\rho \bm{u})=0,\\
			   & \partial_t(\rho \bm{u})+\grad_x\cdot(\rho \bm{u}\otimes \bm{u})
		  	   +\grad_x\cdot\int\frac{\xi\otimes\xi}{\omega^2}g \,d\xi=\rho \BBB(t,x).
			  \end{aligned} \right.
		  \end{equation}
		  In particular, assuming that the couple
                  $(\rho,\bm{u})$ remains smooth\footnote{This is the
                    case at least for short times, and we believe that
                    this can be extended to larger time using the
                    dissipative structure of the right hand side of
                    the equation on $\rho \bm{u}$. See related
                    discussion  in \cite{TadmorTan:2014} for the pressureless system.} and that $\rho$ is nonzero, we have the following equation for the evolution of $\bm{u}$:
		  \begin{equation}
		    \label{eq:u}
		    \left\{ \begin{aligned}
		      & \partial_t \bm{u} +\bm{u}\cdot\grad_x\bm{u} + \frac{1}{\rho}\grad_x\cdot\left(\frac{1}{\omega^2} P\right)=\BBB, \\
		      & \, \\
		  	  & \bm{u}(0,x) = \frac{1}{\rho(0,x)} \int_{\RR^d} f_0(x,v) \,v \, dv,
			  \end{aligned} \right.  
		  \end{equation}
		  where we defined $P$ as a ``pressure'' of $g$, namely
		  \begin{equation*}
		     P=\int_{\RR^d} \xi\otimes\xi \,g(\xi) \, d\xi.
		  \end{equation*}
		  
		  We can now choose the definition of $\omega$. As in section \ref{subHomog}, we want this quantity to be a good indicator of the support of $f$, and for the sake of simplicity we also want its definition to yield a simpler equation for $g$. Using the same arguments, we define $\omega$ as the solution to
		  \begin{equation}
		    \label{eq:omega}
		    \left\{ \begin{aligned}
			    & \partial_t \omega +\bm{u}\cdot\grad_x\omega-\omega \AAA =0, \\
			    & \, \\
			    & \omega(0,x) = \omega_0(x).
			  \end{aligned} \right.  
		  \end{equation}
		  
		  Plugging \eqref{eq:u} and \eqref{eq:omega} in \eqref{eq:gAlmostConservative}, we obtain the general system giving the evolution  of $(g, \bm{u}, \omega)$, namely
		  \begin{equation}
		    \label{eq:g}
		    \left\{ \begin{aligned}
		      & \partial_t g 
		        + \nabla_x \cdot \left [ \left (\bm{u} + \frac{\xi}{\omega}\right ) g\right] 
		        + \nabla_\xi \cdot \left [  \left ( - \xi\cdot\nabla_x \bm{u}
		          + \frac{\xi\cdot\nabla_x \omega}{\omega^2}~\xi + \frac{\omega}{\rho}\grad_x\cdot\left(\frac{1}{\omega^2} P\right)\right ) g \right ] = 0, \\
		     & \partial_t \bm{u} +\bm{u}\cdot\grad_x\bm{u} +  \frac{1}{\rho}\grad_x\cdot\left(\frac{1}{\omega^2} P\right)=\BBB(t,x), \\
		     & \partial_t \omega +\bm{u}\cdot\grad_x\omega-\omega \AAA =0, 
		    \end{aligned} \right.
		  \end{equation}
		  the initial condition for this system being given by
		  \begin{equation}
		    \label{eq:g0}
		    \left\{ \begin{aligned}
		      & \bm{u}(0,x) = \frac{1}{\rho(0,x)} \int_{\RR^d} f_0(x,v) \,v \, dv, \quad \omega(0,x) = \omega_0(x) >0, \\
		      &  g(0,x,\xi) = \frac1{\omega(0,x)^d} f_0\left (\bm{u}(0,x) + \frac\xi{\omega(0,x)}\right ).
		    \end{aligned} \right.
		  \end{equation}
		  Note that the equation for $g$ is in a conservative form, which is of great interest for numerical purposes.
		  
		  \begin{remark}
		    An important feature of the system \eqref{eq:g} is that it is nonlocal only in $x$, through the equation for $\bm u$, whereas the equation \eqref{eq:flocking} for $f$ is nonlocal both in $x$ and $v$. 
		  This provides huge gains in computational time for the new rescaled model.
		  \end{remark}
		  
		  \begin{remark}
		    We can also write the $g$ equation in \eqref{eq:g} as
		    \begin{equation}\label{eq:galt}
		      \partial_t g + \grad_x\cdot(\bm{u}
                      g)-\grad_\xi\cdot\big(\xi\cdot\grad_x\bm{u}~ g\big)+\frac1\omega R=0,
		    \end{equation}
		    where the ``remainder'' term $R$ is given by 
		    \[
		      R=\frac{\omega^2}{\rho}\left(\grad_x\cdot \frac{P}{\omega^2}\right)\cdot\grad_\xi g+{\xi}\cdot\left[\frac{\grad_x\omega} {\omega} \grad_\xi \cdot (\xi g)+\grad_xg\right].
        \]
        Because $\omega$ grows exponentially in time, and $R$ is of
        order 1, the last term on the previous equation on $g$ can be neglected for large time.
		  \end{remark}

			\begin{remark}
			  The coupled system \eqref{eq:g} is in some sense easier to deal with numerically than if one had to use the uncoupled approach introduced in \cite{FilbetRey:2012}. 
			  Indeed, this previous work required the knowledge of a closed macroscopic description of the system. If this is manageable for the Boltzmann equation or for the granular gases equation, here it is more difficult.
			  Indeed, the equilibria of equation \eqref{eq:flocking} being monokinetic
			  \[
			    f(x,v) = \rho(x) \delta(v-\bm{u}(x)), 
			  \] 
			  plugging such a function into equation \eqref{eq:flocking} and computing the zeroth and first moments, one obtain the following dynamics for $\rho$ and $\bm{u}$:
			  \begin{equation}\label{eq:macro}
			  	\left \{\begin{aligned}
				   & \partial_t \rho+\grad_x\cdot(\rho \bm{u})=0,\\
				   & \partial_t(\rho \bm{u})+\grad_x\cdot(\rho \bm{u}\otimes \bm{u})=\BBB.
				  \end{aligned} \right.
                            \end{equation}
				When $\BBB=0$, the equation is usually known as the pressureless Euler system, and exhibits some complicated behavior such as the creation of $\delta$-singularity in finite time \cite{ERykovSinai:1996}.
                                With the alignment force $\BBB$, the
                                solution is less singular. The system has been
                                studied in \cite{TadmorTan:2014},
                                where a critical threshold phenomenon
                                is addressed: subcritical initial data
                                leads to global smooth solution, while
                                supercritical initial data drives to
                                finite time generation of $\delta$-shock.
			\end{remark}

      We end this section by verifying that the zero momentum property \eqref{eq:scalingMomentum} on $g$ is embedded in system \eqref{eq:g0}, which is  an important feature of the equation and will be needed when designing the numerical method.
      \begin{proposition}
        \label{prop:0mom}
        Assume $g=g(t,x,\xi)$ is a smooth solution of \eqref{eq:g} with zero initial momentum 
        \[ 
          M(0,x)=\int\xi \, g(0,x,\xi) \, d\xi=0, \ \forall x\in \Omega.
        \] 
        Then, $M(t,x) =0$ for all $t > 0$ and all $x\in \Omega$.
      \end{proposition}
      
      \begin{proof}
        Multiplying \eqref{eq:galt} by $\xi$ and integrating with
        respect to $\xi$, we have after integration by part that
        \[\partial_tM =
          -\grad_x\cdot(uM)-M\cdot\grad_x\bm{u}-\frac{1}{\omega}\int_{\R^d}\xi \,
          R \,d\xi.\]
        It then suffices to check that 
        \[\int_{\R^d}\xi \, R \,d\xi=0.\]
        Indeed, we have componentwise that
        \begin{align*}
          \int_{\R^d}\xi_kR \,d\xi=
          &\sum_{i,j}\left[\frac{\omega^2}{\rho}\partial_{x_j}\left(
            \frac{P_{ij}}{\omega^2}\right)\int_{\R^d}\xi_k\partial_{\xi_i}
          g\,d\xi+\int_{\R^d}\left(\frac{\partial_{x_j}\omega} {\omega} \partial_{\xi_i}
            (\xi_i g)+\partial_{x_j}g\right)\xi_j\xi_k \,d\xi\right]\\
        =
          &\sum_{i,j}\left[-\frac{\omega^2}{\rho}\partial_{x_j}\left(
            \frac{P_{ij}}{\omega^2}\right)\delta_{ik}\rho
          -\frac{\partial_{x_j}\omega} {\omega}
          (\delta_{ij}P_{ki}+\delta_{ik}P_{ji})+\partial_{x_j}
          P_{kj}\right]\\
        =
          &\sum_{i,j}\left[-\left(\partial_{x_j}P_{kj}-\frac{2\partial_{x_j}\omega}{\omega}P_{kj}\right)
          -\frac{\partial_{x_j}\omega} {\omega}
          (P_{kj}+P_{kj})+\partial_{x_j}
          P_{kj}\right]=0.
        \end{align*}
      \end{proof}

  \section{Numerical Schemes}
    \label{secNumSchemes}
    
    In this section, we present the numerical implementation of the equations for $\bm u$ and $\omega$ in the rescaled dynamics \eqref{eq:g}.

		\subsection{Evolution of the Macroscopic Velocity}
		  \label{subSysConsLawsNum}

	  	We shall solve the dynamics of $\bm{u}$ through the conservative form \eqref{eq:sysMomentCoupled}.
		  More generally, we shall focus on the space discretization of the system of $n$ conservation laws
		  \begin{equation}
			  \label{sysConservLaws}
			  \left \{ \begin{aligned}
			   & \frac{\partial U}{\partial t} + \nabla_x \cdot G(U) = \mathcal H(t,x,U), \ \forall \, (t,x) \in \RR_+ \times \Omega, \\
			   & \, \\
			   & U(0,x) = U_0(x), 
			   \end{aligned} \right.
			\end{equation}
			for a smooth function $G : \RR^n \to M_{n\times d}(\RR)$ and a Lipschitz-continuous domain $\Omega \subset \RR^{d}$. 
			The source term $\mathcal H$ will be problem dependent, and be treated separately.
			Indeed, this equation covers the systems of conservation laws of type \eqref{eq:sysMomentCoupled} with $n=d+1$, $U = (\rho, \rho \bm{u})^\intercal$ and $G$ non linear, or the equation \eqref{eq:omega} ($n=1$, $G$ linear) describing the evolution of $\omega$ for the Cucker-Smale case.		
			Our approach of the problem will be made in the framework of finite volume schemes, using central Lax Friedrichs schemes with slope limiters (see \eg Nessyahu and Tadmor \cite{NessyahuTadmor:1990}).  
			We shall  present the spatial discretization of \eqref{sysConservLaws} in one space dimension for simplicity purposes. The extension for Cartesian grid in the multidimensional case will then be straightforward.

			In the one dimensional setting, the domain $\Omega = (a, b)$ is a finite interval of $\RR$. 
			We define a mesh of $\Omega$, not necessarily uniform, by introducing a sequence of $N_x$ control volume $K_i := \left (x_{i-\frac{1}{2}}, x_{i + \frac 12}\right )$ for $i \in \llbracket 0, N_x \rrbracket$ with $x_i := \left (x_{i-\frac{1}{2}}+ x_{i + \frac 12}\right )/ \,2$ and 
			\begin{equation*}
			  a = x_\frac{1}{2} < x_1 <  \cdots < x_{i-\frac{1}{2}} < x_i < x_{i + \frac 12} < \cdots < x_{N_x} < x_{N_x + \frac 12} = b.
			\end{equation*}
			The Lebesgue measure of the control volume is then simply $\mes(K_i) = x_{i+\frac 12} - x_{i - \frac 12}$. 
			Let $u_i = u_i(t)$ be an approximation of the mean value of $u$ over a control volume $K_i$. By integrating the transport equation \eqref{sysConservLaws} over $K_i$, we get the semi-discrete equations
		  \begin{equation}
			  \label{eqFreeTranspDisc}
			  \left \{ \begin{aligned}
			   & \mes(K_i) \frac{\partial u_i}{\partial t}(t) + {F}_{i+\frac 12} - {F}_{i-\frac 12} = \mathcal{H}_i(t), \ \forall \, t \in \RR_+, \ i \in \llbracket 1, N_x \rrbracket, \\
			   & \, \\
			   & u_i(0) = \frac{1}{\mes(K_i)} \int_{K_i} u_0(x) \, dx.
			   \end{aligned} \right.
			\end{equation}
			
		  In the Cauchy problem \eqref{eqFreeTranspDisc}, the quantity $\left ({F}_{i+\frac 12}\right )_i$, the \emph{numerical flux}, is an approximation of the flux function $x \mapsto G\left (u(t,x)\right )$ at the cell interface $x_{i+\frac 12}$. We choose to use the so-called Lax-Friedrichs fluxes with the second order Van Leer's slope limiter \cite{VanLeer:1977a}. 
		  In this setting, the slope limited  flux is given by 
			\begin{equation*}
			  {F}_{i+\frac 12} = \frac12 \left ( G\left (u_{i + \frac 12, -}\right ) + G\left (u_{i + \frac 12, +}\right )\right ) - \frac{\lambda_{i+\frac 12}}{2}\left ( u_{i + \frac 12, -} - u_{i + \frac 12, +}\right )
			\end{equation*}
			where we have set
			\[
			  \lambda_{i+\frac 12} := \max_{\lambda \in {\rm Sp \left (G'(U_i)\right )}} |\lambda|, 
			\]
			and $u_{i + \frac 12, \pm}$ is the slope limited reconstruction of $u$ at the cell interface, namely, componentwise,
			\begin{equation*}
				\left \{ \begin{aligned}
				  u_{i + \frac 12, -} & = u_i + \frac 12 \, \phi( \theta_i ) \, ( u_{i+1} - u_i), \\
				   & \, \\
				  u_{i + \frac 12, +} & = u_i - \frac 12 \, \phi( \theta_i ) \, ( u_{i+2} - u_{i+1}).
				\end{aligned} \right.
			\end{equation*}
			In this last expression, $\theta_i$ is the \emph{slope} for each component $u_k$ of $u$:
			\[ \theta_{i,k} = \frac{u_{i,k} - u_{i-1,k}}{u_{i+1,k} - u_{i,k}}, \quad \forall k \in \{1, \ldots,n\} \]
			and $\phi$ is the so-called \emph{Van Leer}'s limiter
			\[ \phi(\theta) := \frac{\theta + |\theta|}{1 + \theta}. \]
			
			We notice that this second order method uses a $2$ points stencil, and we will then have to define the value of the solution on the \emph{ghost cells} 
			\[ \left \{x_{-\frac32},x_{-\frac12},x_{N_x+\frac32},x_{\frac52}\right \}.\]
			This value will be set according to the boundary conditions chosen for the problem at hand.

	  \subsection{Discretization of the Flocking Terms}
	    \label{subDiscFlockingTerm}
	    
	    Since we are only dealing with first order schemes for the transport parts in \eqref{eq:g}, we shall not use a high order spectral method for the discretization of the flocking terms $\AAA$ and $\BBB$.
	    We then simply approximate these terms using a first order quadrature rule. On an uniform grid $x_i = i \Delta x$ and $\xi_j = j \Delta \xi$ for $\Delta x >0$ and $\Delta \xi$, we have for the Cucker-Smale model \ref{eq:AB_CS}:
	    \[
	      \AAA_{CS,i} = \Delta x \,\sum_j \phi\left (|x_i - x_j|\right ) \rho_i, \quad \BBB_{CS,i} = \Delta x \,\sum_j \phi\left (|x_i - x_j|\right ) \left ( \bm{u}_j - \bm{u}_i\right ) \rho_i,
	    \]			
	    where we set
	    \[ \rho_i = \Delta \xi \sum_j g_{ij}, \quad \rho_i \bm{u}_i = \Delta \xi \sum_j \xi_j \, g_{ij}.\]
	    The Motsch-Tadmor model is then simply given by
	    \[
	    \AAA_{MT,i} = 1, \quad \BBB_{MT,i} = \frac{\BBB_{CS,i}}{\AAA_{CS,i}}.
	    \]

\subsection{Evolution of the scaling factor}
Let us recall the dynamics of the scaling factor $\omega$
\[\partial_t\omega+\bm{u}\cdot\grad_x\omega-\omega\AAA=0.\]

In the Motsch-Tadmor setup, we have seen that 
\[\AAA_{MT}\equiv1.\] 
If one pick the initial
scaling $\omega_0=1$, then there is an explicit spatially homogeneous
solution for $\omega$, given by:
\[\omega(x,t)=e^t.\]

In the Cucker-Smale setup, $\omega$ is spatially dependent. Along the
characteristic flow,
\[\omega'=\omega\AAA_{CS},\qquad \text{where } '=\partial_t+\bm{u}\cdot\grad_x.\]
If $\AAA_{CS}=\phi\star\rho$ is strictly positive, $\omega$ grows exponentially in time.
From theorem \ref{thm:flocking}, we get a uniform in space-time lower bound on
$\AAA_{CS}$:
\[\AAA_{CS}(x,t)\geq\phi(D)\|\rho\|_{L^1},\]
where $D=\sup_tS(t)$ is finite and $\phi(D)>0$ in the case of
flocking. Hence, $\omega$ has an exponential growth as well for
Cucker-Smale system.

To evolve $\omega$ numerically, we rewrite the equation in the conservative form
\[\partial_t\omega+\grad_x\cdot(\bm{u}\omega)=
\omega(\grad_x\cdot\bm{u}+\AAA_{CS}).\]
and it can be treated under the framework of system
\eqref{sysConservLaws}.

	\section{A Momentum Preserving Correction of the Upwind Scheme}
	  \label{secNewFluxes}
	  
	  We have seen in Proposition \eqref{prop:0mom} that one of the important features of the rescaled equation \eqref{eq:g} is that it preserves the zero momentum condition of $g$:
	  \[
	    M(t,x) := \int_{\RR^d} \xi \, g(t,x,\xi) \, d\xi = 0, \ \forall t > 0 \quad \text{if} \quad M(0,x) = 0, \ \forall x \in \Omega.
	  \]
	  We will derive in this section a numerical method that is able to propagates \emph{exactly} this particular property, for the type of equation we are dealing with\footnote{An extension to a more general class of equations is currently in progress \cite{ReyTan:2014b}.}.

		\subsection{A Toy Model}
		  \label{sub:ToyModel}
		  
		  We will start by presenting our approach on a toy model. Let us consider, for $c \in \RR$ constant, the following transport equation on $g=g(t,\xi)$:
			\begin{equation}\label{eq:toy}
			  \left \{ \begin{aligned}
			    & \partial_tg+c \, \grad_\xi\cdot(\xi g)=0,\\
			    & \, \\
			    & g(0,\xi)=g_0(\xi) > 0,
			  \end{aligned} \right.
			\end{equation}
			with the zero initial momentum property 
			\[ \int_{\RR^d}\xi \, g_0(\xi) \,d\xi=0. \]
			A simple calculation yields that if $g$ is solution to \eqref{eq:toy} then one has
			\[\frac{d}{dt}\int_{\RR^d} \xi \,g(t,\xi)\,d\xi=-c\int_{\RR^d}\grad_\xi\cdot(\xi
			g)\,d\xi=c\int_{\RR^d}\xi \,g(t,\xi)\,d\xi.\]
			Therefore, momentum is conserved in time:
			\begin{equation}
			  \label{eq:evolMomentumExactToy}
			  M(t) = e^{ct}\int_{\RR^d}\xi \,g_0 \,d\xi=0.
			\end{equation}
			
			Let us now consider the numerical approximation of this problem. 
			One classical way to consider it is to apply the classical upwind scheme 
			(denoted in all the following by \upwind) to solve the equation. 
			For the sake of simplicity, let us take $c=1$ and a one dimensional, equally distributed grid on $\xi$:
			\[ \xi_j=(j-J)\Delta\xi, \quad \forall j \in \ZZ, \] 
			for a given $J \in \RR$. In particular, 0 is on the grid:
		  \[\xi_J=0.\] 
		  In this case, the fully discrete upwind scheme reads \cite{leveque:2002}
			\begin{equation}\label{eq:FV}
			  g^{n+1}_j=g^n_j-\frac{\Delta t}{\Delta\xi}\left(F^n_{j+1/2}-F^n_{j-1/2}\right),
			\end{equation}
			where the numerical flux $\left (F^n_{j+1/2}\right )_j$ is given since $g^n_j$ is nonnegative by
			\begin{equation}
				\label{eq:classicalFlux}
				F^n_{j+1/2}=\begin{cases}
				\xi_{j+1/2}\,g_j^n& j\geq J\\
				\xi_{j+1/2}\,g_{j+1}^n&j\leq J-1
				\end{cases}.
			\end{equation}
		
			Let us compute the evolution of the \emph{discrete} momentum $M^n$:
			\[M^n := \Delta \xi\sum_j\xi_j\,g_j^n \simeq M\left (t^n\right ).\]
			Since $\xi_J = 0$, the contribution of the flux can be simplified as follows.
			\begin{align*}
				\sum_j\xi_jF^n_{j+1/2}=&
				\sum_{j\geq  J}\xi_j\xi_{j+1/2}\,g_j^n+
				\sum_{j\leq  J-1}\xi_j\xi_{j+1/2} \,g_{j+1}^n\\
				=&\sum_{j\geq J}\xi_j\left(\xi_j+\frac{\Delta\xi}{2}\right)g_j^n+
				\sum_{j\leq J}(\xi_j-\Delta\xi)\left(\xi_j-\frac{\Delta\xi}{2}\right)g_j^n\\
				=&\sum_j\xi_j^2 \,g_j^n+\frac{\Delta\xi}{2}\left[\sum_{j\geq  J}\xi_j\,g_j^n
				-3\sum_{j\leq J}\xi_j\,g_j^n\right]+\frac{(\Delta\xi)^2}{2}\sum_{j\leq J}g_j^n.
			\end{align*}
			Similarly, we get
			\[
				\sum_j\xi_jF^n_{j-1/2}=\sum_j\xi_j^2\, g_j^n+\frac{\Delta\xi}{2}\left[3\sum_{j\geq  J}\xi_j\,g_j^n
				-\sum_{j\leq J}\xi_j\,g_j^n\right]+\frac{(\Delta\xi)^2}{2}\sum_{j\geq J}g_j^n.
			\]		
			The evolution of the discrete momentum then reads,
			\begin{align}
				\sum_j\xi_j\,g^{n+1}_j=&\sum_j\xi_j\,g^n_j-\frac{\Delta
				  t}{\Delta\xi}\left[\sum_j\xi_jF^n_{j+1/2}-\sum_j\xi_jF^n_{j-1/2}\right] \notag \\
				=&(1+\Delta t)\sum_j\xi_j\,g^n_j+\frac{\Delta t\Delta\xi}{2}
				\left[-\sum_{j\leq J} g_j^n+\sum_{j\geq J}g_j^n\right], \label{eq:evolMomentumUpwind}
			\end{align}
			namely one has
			\[ 
			  M^{n+1} = (1 + \Delta t) M^n + \frac{\Delta t\Delta\xi}{2}
			  \left[-\Delta\xi\sum_{j\leq J} g_j^n+\Delta\xi\sum_{j\geq J}g_j^n\right].
			\]		
			In the special case where $g$ is symmetric in $\xi$, the discrete zero
			momentum is preserved in time. However, it is in general not true, unless one has
			\[\sum_{j\leq J} g_j^n=\sum_{j\geq J}g_j^n.\]
		
			To ensure momentum conservation, we introduce a correction $\tilde{F}^n$ on the
			upwind flux $F^n$. If the correction satisfies
			\begin{equation}\label{eq:corr}
		  	\sum_j\xi_j(\tilde{F}_{j+1/2}^n-\tilde{F}_{j-1/2}^n)=
		  	\frac{(\Delta\xi)^2}{2}\left[\sum_{j\leq J} g_j^n-\sum_{j\geq
			      J}g_j^n\right],
			\end{equation}
			then the new flux $F^n+\tilde{F}^n$ will preserve zero momentum, according to \eqref{eq:evolMomentumUpwind}.	
			
      We provide two corrections $\tilde{F}^{(1),n}$ and $\tilde{F}^{(0),n}$
		  which satisfy \eqref{eq:corr}:
			\[
			  \tilde{F}^{(1),n}_{j+1/2}=\begin{cases}
			    -\frac{\Delta\xi}{2}g^n_{j+1}& j\geq J\\
			    +\frac{\Delta\xi}{2}g^n_j&j\leq J-1
			  \end{cases},\qquad
	    	\tilde{F}^{(0),n}_{j+1/2}=\begin{cases}
		    	-\frac{\Delta\xi}{2}g^n_j& j\geq J\\
		    	+\frac{\Delta\xi}{2}g^n_{j+1}&j\leq J-1
		  	\end{cases}.
		  \]
      The corresponding new fluxes $F^{(1),n}$ and $F^{(0),n}$ have the following forms.
			\begin{equation}\label{eq:newflux}
		  	F^{(1),n}_{j+1/2}=\begin{cases}
		    	\xi_{j+1/2}\,g^n_j-\frac{\Delta\xi}{2}g^n_{j+1}& j\geq J\\
		    	\xi_{j+1/2}\,g^n_{j+1}+\frac{\Delta\xi}{2}g^n_j&j\leq J-1
		  	\end{cases},\qquad
	      F^{(0),n}_{j+1/2}=\begin{cases}
		    	\xi_j\,g^n_j& j\geq J\\
		    	\xi_{j+1\,}g^n_{j+1}&j\leq J-1
		  	\end{cases}.
	    \end{equation}
			Moreover, we can get a family of corrections satisfying \eqref{eq:corr} by
			interpolating between $\tilde{F}^{(1),n}$ and $\tilde{F}^{(0),n}$:
			\[
			  \tilde{F}^{(\theta),n} :=\theta\tilde{F}^{(1),n}
		  	+(1-\theta)\tilde{F}^{(0),n},\quad \theta\in[-1,1].
		  \]
			The respective fluxes $F^{(\theta),n}:=F^n+\tilde{F}^{(\theta),n}$ can be then expressed as
			\begin{equation}\label{eq:newfluxes}
		  	F^{(\theta),n}_{j+1/2}=F^{(0),n}_{j+1/2}-\frac{\theta\Delta\xi}{2}(g_{j+1}^n-g_j^n).
		  \end{equation}
		  We have the following result:
		
			\begin{proposition}\label{prop:conserve}
			  Suppose that the sequence $\{g^0_j\}_j$ has $0$ discrete momentum:
			  \[\Delta \xi \sum_j\xi_j\,g^0_j=0.\] 
			  Then for any $\theta \in [-1,1]$, the scheme \eqref{eq:FV} with initial condition $\{g^0_j\}_j$ 
			  with numerical fluxes $F^{(\theta)}$ given by \eqref{eq:newfluxes} preserves 
			  the discrete mass $\Delta \xi \sum_j g_j^n$ and momentum $\Delta \xi \sum_j\xi_j\,g_j^n$.
			\end{proposition}
		
			For the general case $c\in\RR$, a similar correction can be added into 
			the upwind flux:
			\begin{equation}\label{eq:newfluxgeneral}
	      F^{(\theta),n}_{j+1/2}=F^{(0),n}_{j+1/2}-\frac{c\theta\Delta\xi}{2}(g_{j+1}-g_j),
			\end{equation}
			where $F^{(0)}$ is defined as
			\[
			  F^{(0),n}_{j+1/2}=\begin{cases}
			    c\,\xi_j\,g^n_j& j\geq J\\
			    c\,\xi_{j+1}\,g^n_{j+1}&j\leq J-1
			  \end{cases} \text{for } c>0,\quad\text{and}
			  ~~
			  F^{(0),n}_{j+1/2}=\begin{cases}
			    c\,\xi_j\,g^n_j& j\leq J-1\\
			    c\,\xi_{j+1}\,g^n_{j+1}&j\geq J
			  \end{cases}\text{for } c<0.
			\]
			In the following, this family of fluxes will be called 
			\mcu{$\theta$}, for \emph{Momentum Conservative Upwind}
			fluxes.
			With the correction, it is easy to check that
			\[M^{n+1}=(1+c\Delta t)M^n.\]
			It implies that the zero momentum property is preserved. Moreover, even if the
			initial momentum is not zero, the new flux provides a good
			approximation of the momentum. 
			Indeed, note that $1+c\Delta t$ is a
			first-order in time approximation of $e^{c\Delta t}$, the correct 
			behavior of the momentum of a solution to \eqref{eq:toy}, 
			according to \eqref{eq:evolMomentumExactToy}. It is moreover
			independent with the choice of $\Delta\xi$. Various ways can be
			applied to obtain higher time accuracy.
			
      Let us summarize the properties of the \mcu{$\theta$} fluxes.
			\begin{proposition}
			  Consider a finite volume scheme \eqref{eq:FV} for the approximation of equation \eqref{eq:toy}
			  with the numerical fluxes \mcu{$\theta$}, for a given $\theta\in[-1,1]$. Then one has:
				\begin{enumerate}
		  		\item \emph{Accuracy.} The scheme solves the equation
				    \eqref{eq:toy} with first order accuracy.
		  		\item \emph{Positivity preserving.} If $c\,\theta\geq0$ and the computational domain is
			  	  $[-L,L]$, then the scheme preserves positivity
			  	  under the CFL condition
			    	\begin{equation}
			    	  \label{eq:CFLcond}
			    	  \lambda=\frac{\Delta t}{\Delta\xi}\leq\frac{1}{|c|L}.
			    	\end{equation}
			  	\item \emph{Mass conservation.} The scheme preserves the discrete mass $\Delta \xi \sum_j g_j^n$.
			  	\item \emph{Momentum conservation.} The scheme preserves the zero initial
			  	  momentum  property.
				\end{enumerate}
			\end{proposition}

			\begin{remark}
				As a direct consequence of positivity preserving and mass
				conservation, the scheme is $l^1$-stable, namely, the discrete $l^1$
				norm $\|g_j\|_{l^1}$ is conserved in time if $c \, \theta \geq 0$. In particular, \mcu{0} is
				stable for any choice of $c$, under the CFL condition \eqref{eq:CFLcond}.
			\end{remark}
		
			\begin{remark}
			  The new flux can be easily implemented in a higher dimensional setting, as in each
		  	interface, the fluxes can be treated like in the $d=1$ case.
			\end{remark}
%

		\subsection{Application to Flocking Models}
		  \label{sub:flockMod}

      \paragraph{\emph{The Motsch-Tadmor Dynamics.}}
			Let us now present the discretization of the equation describing the evolution of $g$.
			We shall apply the new momentum preserving flux to solve this equation for different
			models, starting with Motsch-Tadmor. Here, we recall the dynamics in 1D
			\begin{equation}\label{eq:gMT1D}
			  \partial_t g +\partial_x \left [ \left (u + \frac{\xi}{\omega}\right ) g\right] 
					        + \partial_\xi \left [  \left ( - \xi\,\partial_xu
					        + \frac{1}{\rho\omega}\partial_x P\right ) g
			                    \right ] = 0.
			\end{equation}
			We consider $x\in\mathbb{T}$ and an initial density bounded by below:
			\[ \int_{\RR^d} f_0(x,v) \, dv \geq \rho_0 >0.\] 
			Thus, because of the continuity equation, vacuum cannot exist in
			any finite time if $u_x$ remains bounded.  

			We treat the four terms arising in \eqref{eq:gMT1D} one by one using a
			finite volume method:
			\begin{align*}
				g_{ij}^{n+1}=g_{ij}^n&-\frac{\Delta t}{\Delta x}\left(F^{1,n}_{i+1/2,j}-F^{1,n}_{i-1/2,j}\right)
				-\frac{\Delta t}{\Delta x}\left(F^{2,n}_{i+1/2,j}-F^{2,n}_{i-1/2,j}\right)\\
				&-\frac{\Delta t}{\Delta\xi}\left(F^{3,n}_{i,j+1/2}-F^{4,n}_{i,j-1/2}\right)
				-\frac{\Delta t}{\Delta\xi}\left(F^{4,n}_{i,j+1/2}-F^{4,n}_{i,j-1/2}\right),
			\end{align*}
			where $F^1,\ldots,F^4$ are the numerical fluxes associated respectively to
			\[ 
			  x \mapsto u\,g, \quad x\mapsto\frac{\xi}{\omega}\,g, \quad 
			  \xi \mapsto -\xi\,\partial_xu\,g, \quad\xi \mapsto \frac{1}{\rho\omega}\partial_x
			P g.
			\] 
			For the sake of simplicity, we discretize $g$ with equally
			distributed cells in both $x$ and $\xi$. In particular, $\xi_J=0$.
			We shall also omit the superscript $n$ from now on.

			For $F^1$, we take the classical \upwind~flux \eqref{eq:classicalFlux}
			\[ F^1_{i+1/2,j}=\begin{cases}
	  		u_{i+1/2}\,g_{ij}&u_{i+1/2}\geq0\\
  			u_{i+1/2}\,g_{i+1,j}&u_{i+1/2}<0
			\end{cases},\]
			which clearly preserves the discrete momentum.
 
			Similarly, the classical \upwind~flux can be used for $F^2$ and $F^4$ as
			well,
			\[F^2_{i+1/2,j}=\begin{cases}
				\displaystyle\frac{\xi_j}{\omega}\,g_{ij}&j\geq J\\
				\, \\
				\displaystyle\frac{\xi_j}{\omega}\,g_{i+1,j}&j\leq J-1
				\end{cases},\qquad F^4_{i,j+1/2}=\begin{cases}
				\displaystyle\frac{1}{\rho_i\omega}(\partial_xP)_i\,g_{ij}&(\partial_xP)_i\geq0\\
				\, \\
				\displaystyle\frac{1}{\rho_i\omega}(\partial_xP)_i\,g_{i,j+1}&(\partial_xP)_i<0
			\end{cases}.\]
			As these two terms preserves momentum when combined together, namely,
			\[\int_\R\xi\left[\partial_x\left(\frac{\xi}{\omega}\,g\right)
	  		+\partial_\xi\left(\frac{1}{\rho\omega}\partial_xPg\right)\right]d\xi=0,\]
			the discrete momentum will be preserved with the discrete flux as well:
			\[\sum_j\xi_j\left[\frac{F^2_{i+1/2,j}-F^2_{i-1/2,j}}{\Delta x}+
	  		\frac{F^4_{i,j+1/2}-F^4_{i,j-1/2}}{\Delta\xi}\right]=0.\]
			This can be easily checked if we approximate $\partial_xP$ and $\rho$
			by
			\[(\partial_xP)_i=\Delta\xi\left[
			  \sum_{j\geq
                            J+1}\xi_j^2\frac{g_{ij}-g_{i-1,j}}{\Delta
                            x}+\sum_{j\leq
                            J-1}\xi_j^2\frac{g_{i+1,j}-g_{ij}}{\Delta x}
		  	\right],\quad\rho_i=\Delta\xi\sum_jg_{ij}.\]
			
			Finally, the last remaining term $F^3$ can be reduced to the
			toy model for fixed $i$, where the coefficient $c=-(\partial_xu)_i$
			is $x$-dependent. 
			Hence, we apply the new \mcu{$\theta$} flux \eqref{eq:newfluxgeneral} 
			and the zero momentum property is
			preserved as a consequence of proposition \ref{prop:conserve}. 
			Moreover, concerning the question of
			stability, we have seen that \mcub~preserves the positivity for all $c \in \RR$. 
			Alternatively, we can also use \mcua~for $c>0$ and \mcu{-1} for $c<0$.

                        \medskip
      \paragraph{\emph{The Cucker-Smale Dynamics.}}
      We now consider the Cucker-Smale model. In this case, the scaling factor
      $\omega$ depends on both time and space variables, which brings an extra term
      to the $g$ equation, as well as some new numerical difficulties. 
      The 1D dynamics reads
			\begin{equation}\label{eq:gCS1D}
			  \partial_t g +\partial_x \left [ \left (u + \frac{\xi}{\omega}\right ) g\right] 
					        + \partial_\xi \left [  \left ( -\xi\,\partial_xu
					        + \frac{\omega}{\rho}\partial_x\left(\frac{P}{\omega^2}\right)
			                        + \frac{\partial_x\omega}{\omega^2}\xi^2\right ) g
			                    \right ] = 0.
			\end{equation}

			For its numerical discretization, the quantities $F^1$ and $F^3$ are the same as
			in the Motsch-Tadmor case. The quantity $F^2$ is again treated by the upwind flux
      \[F^2_{i+1/2,j}=\begin{cases}
				\displaystyle\frac{\xi_j}{\omega_{i+1/2}}g_{ij}&j\geq J\\
				\displaystyle\frac{\xi_j}{\omega_{1+1/2}}g_{i+1,j}&j\leq J-1
			\end{cases}.\]
			To get $\omega_{i+1/2}$, we can either evolve $\omega$ on a staggered grid
			$\{x_{i+1/2}\}_i$, or interpolate from the knowledge of the cell-centered values
			$\{\omega_i\}_i$. We choose in all our numerical experiments this latter approach,
			with a simple first order interpolation.

			Since $\omega$ is not a constant, we also have to modify $F^4$
			as follows:
			\[F^4_{i,j+1/2}=\begin{cases}
				\displaystyle\frac{1}{\rho_i}\frac{\omega_{i+1/2}+\omega_{i-1/2}}{2}\left[\partial_x\left(\frac{P}{\omega^2}\right)\right]_ig_{ij}
				&\displaystyle\left[\partial_x\left(\frac{P}{\omega^2}\right)\right]_i\geq0\\
				&\\
				\displaystyle\frac{1}{\rho_i}\frac{\omega_{i+1/2}+\omega_{i-1/2}}{2}\left[\partial_x\left(\frac{P}{\omega^2}\right)\right]_ig_{i,j+1}
				&\displaystyle\left[\partial_x\left(\frac{P}{\omega^2}\right)\right]_i<0
			\end{cases},\]
			where
			\[\left[\partial_x\left(\frac{P}{\omega^2}\right)\right]_i
			  =\frac{\Delta\xi}{\Delta x}\left[
	  		\sum_{j\geq J+1}\xi_j^2\left(\frac{g_{ij}}{\omega_i^2}
	  		-\frac{g_{i-1,j}}{\omega_{i-1}^2}\right)
	  		+\sum_{j\leq J-1}\xi_j^2\left(\frac{g_{i+1,j}}{\omega_{i+1}}
		  	-\frac{g_{ij}}{\omega_i}\right)
			\right],\quad\rho_i=\Delta\xi\sum_jg_{ij}.\]

      Finally, for the additional term 
      \[ \xi \mapsto \frac{\partial_x\omega}{\omega^2}\xi^2g,\] 
      we choose the corresponding flux $F^5$ which is compatible with $F^4$, so
      that discrete zero momentum is preserved. One simple momentum
			preserving flux reads
			\[F_{i,j+1/2}^5=\begin{cases}
		  	\displaystyle\frac{\xi_j^2}{2}\left[\frac{g_{i,j}}{\omega_{i+1/2}^2}+\frac{g_{i-1,j}}{\omega_{i-1/2}^2}\right]\cdot\frac{\omega_i-\omega_{i-1}}{\Delta  x}
		  	& j\geq J\\
	  	  &\\
		  	\displaystyle	\frac{\xi_j^2}{2}\left[\frac{g_{i+1,j}}{\omega_{i+1/2}^2}+\frac{g_{i,j}}{\omega_{i-1/2}^2}\right]\cdot\frac{\omega_{i+1}-\omega_i}{\Delta  x}
			& j\leq J-1\\
			\end{cases}.\]
			Despite the preservation of momentum, this flux uses the information in
			the cell $j$ to determine the flux at the interface $j+1/2$, which is not
			promising. Here, we provide a more reasonable
                        momentum preserving flux, namely
			\[F_{i,j+1/2}^5=\begin{cases}
				\displaystyle\frac{\xi_j\xi_{j+1}}{4}\left(\frac{g_{i,j}+g_{i,j+1}}{\omega_{i+1/2}^2}+\frac{g_{i-1,j}+g_{i-1,j+1}}{\omega_{i-1/2}^2}\right)
				\cdot\frac{\omega_i-\omega_{i-1}}{\Delta  x}
				& j\geq J\\
				&\\
				\displaystyle\frac{\xi_j\xi_{j+1}}{4}\left(\frac{g_{i+1,j}+g_{i+1,j+1}}{\omega_{i+1/2}^2}+\frac{g_{i,j}+g_{i,j+1}}{\omega_{i-1/2}^2}\right)
				\cdot\frac{\omega_{i+1}-\omega_i}{\Delta  x}
				& j\leq J-1\\
			\end{cases}.\]

  \section{Numerical Simulation}
    \label{secNumSim}
    
    \subsection{Test 1 - The Anti-Drift Equation}
      \label{subToyMod}
      
      Before presenting numerical simulations for the full flocking
      equation \eqref{eq:flocking}, we will first demonstrate the
      efficiency of the new conservative fluxes \mcu{$\theta$} described in section \ref{secNewFluxes}. 
      For this, we consider the toy model \eqref{eq:toy}, with $c=1$ (also know as the linear \emph{anti-drift} equation):
      \begin{equation}
        \label{eq:antiDrift}
        \partial_t g + \nabla_x \cdot (\xi g) = 0,
      \end{equation}
      with homogeneous Dirichlet boundary conditions.
      We consider the case $d=1$ and take as an initial condition a sum of two Gaussian functions, with $0$ momentum:
      \[
        g(0,\xi) = \frac1{\sqrt{2 \pi T}} \left[ \frac12 \exp\left ( \frac{|\xi - c_1|^2}{2 T}\right ) + \frac32 \exp\left ( \frac{|\xi - c_2|^2}{2 T}\right ) \right ], \quad T = 0.01, \ c_1 = 0.9375, \ c_2 = 0.3125.
      \]
      One can check that this function has $0$ momentum, but is not symmetric.
      We aim to compare the new momentum conservative upwind fluxes with the classical ones.
      
      \begin{figure}
        \includegraphics[scale=1]{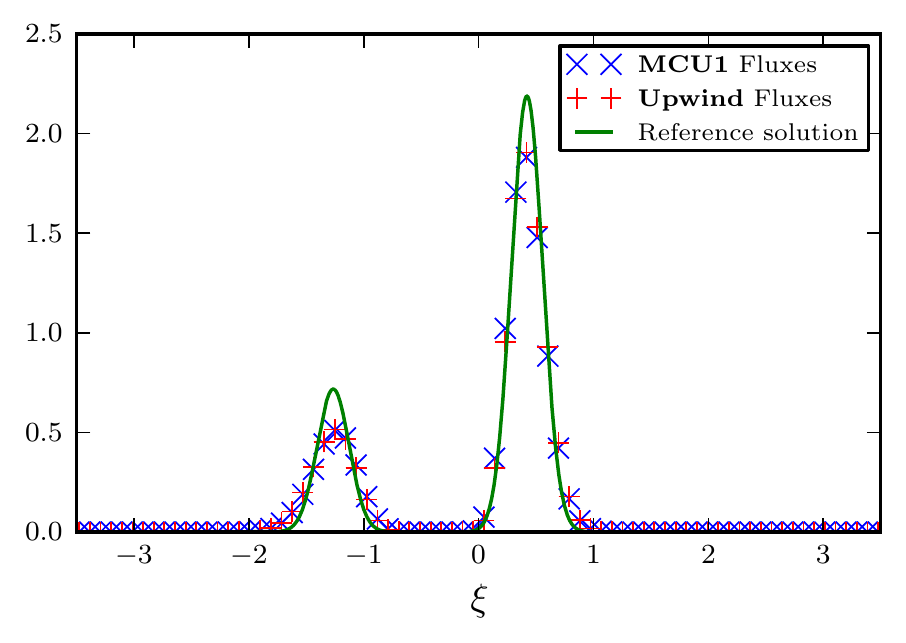}
        \caption{\textbf{Test 1 -} Approximate solutions to the anti-drift equation \eqref{eq:toy} ($c=1$) given by \upwind~and \mcua, at time $t = 0.3$.}
        \label{fig:compSolutions}
      \end{figure}

      We first present in Figure \ref{fig:compSolutions} the approximate solution at time $t=0.3$ of equation \eqref{eq:antiDrift}, obtained with the \upwind~and \mcua~first order fluxes with $N_\xi = 101$ points in the $\xi$ variable to discretize the box $[-3.5,3.5]$. 
      The time stepping is done using a forward Euler discretization with $\Delta t = 1/300$. 
      We also show a reference solution obtained by using a second
      order flux limited scheme as presented in section
      \ref{subSysConsLawsNum} with $2000$ points in $\xi$ and $\Delta t = 1/1500$.
      We observe that both \upwind~and \mcua~fluxes give very similar results, which are in good agreement with the reference solution. Being first order, both schemes are quite diffusive but seem to give the correct wave propagation speed.

      \begin{table}
	      \begin{tabular}{c|ccc}
	        \hline $N_\xi$ & \upwind & \mcua   & \mcub    \\
	        \hline 101     & 4.6e-3  & 3.6e-16 & 3.1e-16  \\
	               201     & 2.3e-3  & 2.1e-16 & 1.9e-16  \\
	               401     & 9.9e-4  & 2.8e-16 & 2.3e-16  \\
	        \hline 
	      \end{tabular}
	      \caption{\textbf{Test 1 -} $L^\infty$ norm of the first moment $M(t)$ of the solution to the drift equation \eqref{eq:toy}, for $t \in [0,0.2]$.}
	      \label{tab:firstMom}
      \end{table}
            
\begin{figure}
        \includegraphics[scale=1]{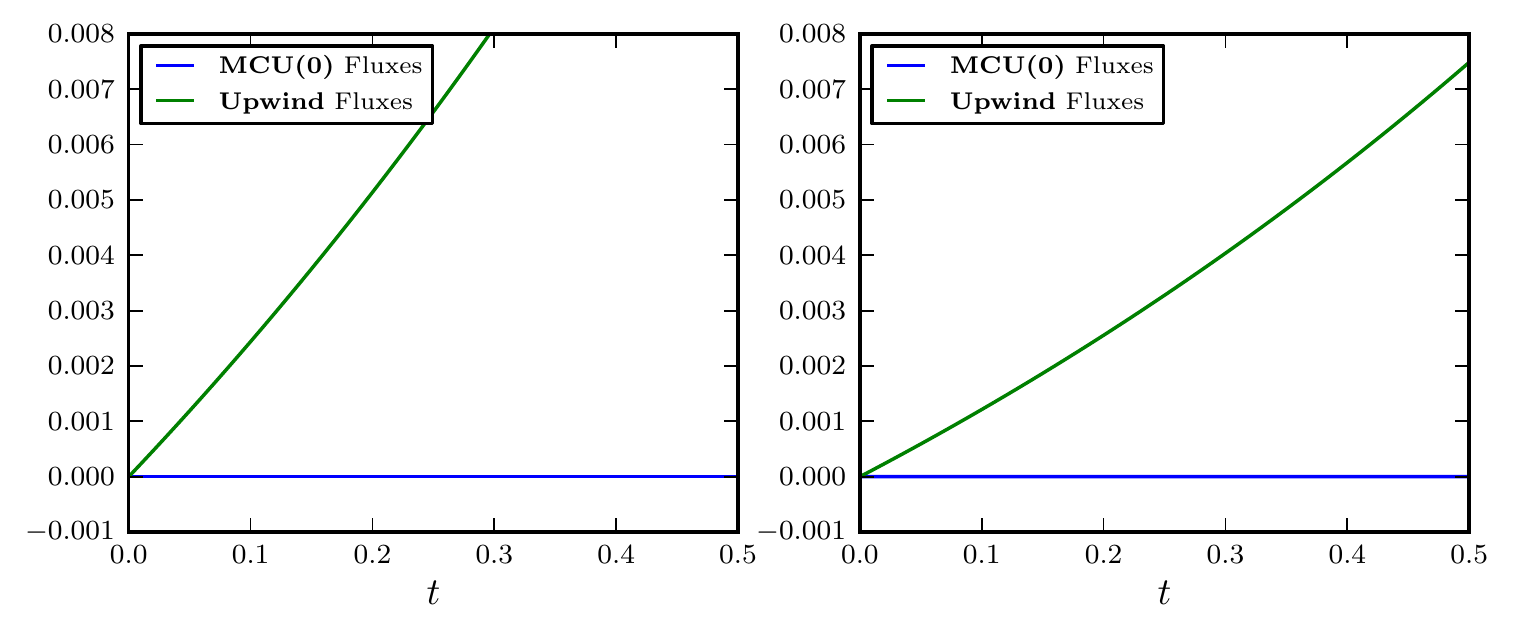}
        \caption{\textbf{Test 1 -} First moment $M(t)$ of the anti-drift equation \eqref{eq:toy} ($c=1$) given by \upwind~and \mcua.
        Coarse grid $N_\xi = 101$ (left) and fine grid $N_\xi = 201$ (right).}
        \label{fig:consVSnoncons}
      \end{figure}

      We then investigate the desired properties, namely the preservation of the first moment of $g$:
      \[ 
        M(t) := \int_{\RR^d} \xi \, g(t, \xi) \, d\xi = 0, \ \forall t \geq 0.
      \]
      We present in Table \ref{tab:firstMom} the $L^\infty$ norm of this quantity for $t \in [0, 0.2]$, for the \upwind, \mcua, and \mcub~fluxes, and for different mesh sizes. 
      We observe that both \mcua~and \mcub~preserves exactly (up to the machine precision) the first moment of $g$, without any influence of the grid size. 
      This is not the case for the classical \upwind~fluxes, where the value of $M(t)$ decreases almost linearly with the size of the mesh.
      This is even more clear in Figure \ref{fig:consVSnoncons}, where we compare the time evolution of the  approximate value of $M(t)$ obtained with the \upwind~and \mcua~fluxes. 
      We take successively $N_\xi = 101$ and $N_\xi = 201$ grid points and respectively $\Delta t = 1/300$ and $\Delta t = 1/600$. 
      While the momentum obtained with the \mcua~fluxes remains nicely $0$ during time, the one obtained with the \upwind~fluxes grows linearly with time.
      Moreover, as expected through equation \eqref{eq:evolMomentumUpwind}, the growth rate of this quantity is proportional to the mesh size.
            
    \subsection{Test 2 - One Dimensional Motsch-Tadmor}
      \label{subNumFlock}
      
      We are now interested in numerical simulations of the flocking equation \eqref{eq:flock} in the Motsch-Tadmor case \eqref{eq:collisionMT} for $d=1$ with the local influence function
      \[
        \phi(r) = \bm 1_{\{|r| \leq 0.1\}},
      \]
      and periodic boundary conditions in $x$.
      We will use for this the rescaled model \eqref{eq:g}. 
      We recall that in the particular Motsch-Tadmor case, the equation describing the evolution for $g$ is reduced to \eqref{eq:gMT1D}. 
      Moreover, we can chose $\omega(t) = \exp(t)$ for all $t \geq 0$ according to \eqref{eq:omegaMT}.
      
      We take as an initial condition the Gaussian function
      \[
        f_0(x,v) = \frac{\rho(x)}{\sqrt{2 \pi}} \exp \left (-|v - u(x)|^2/2\right ), \ \forall x \in \TT, \ v \in \RR,
      \]
      where the density $\rho$ is almost localized in space
      \[
        \rho(x) = 0.01 + \frac{1}{\sqrt{2 \pi T} }\exp(-|x|^2/2T), \quad T = 0.01,
      \] 
      and the momentum $u$ is an oscillating function
      \[
        u(x) = 5 + \sin(x/2 \pi).
      \]
      
      The discretization of the drift part of \eqref{eq:gMT1D} is dealt with using the momentum preserving fluxes \mcub, as presented in section \ref{sub:flockMod}. 
      We take $N_x = 75$ points in the physical space and $N_\xi = 101$ points in the rescaled velocity space for a rescaled velocity variable $\xi \in [-15,15]$. We choose $\Delta t = 1/1500$ because of the size of this support.

      \begin{figure}
        \begin{tabular}{c}
          \includegraphics[scale=1]{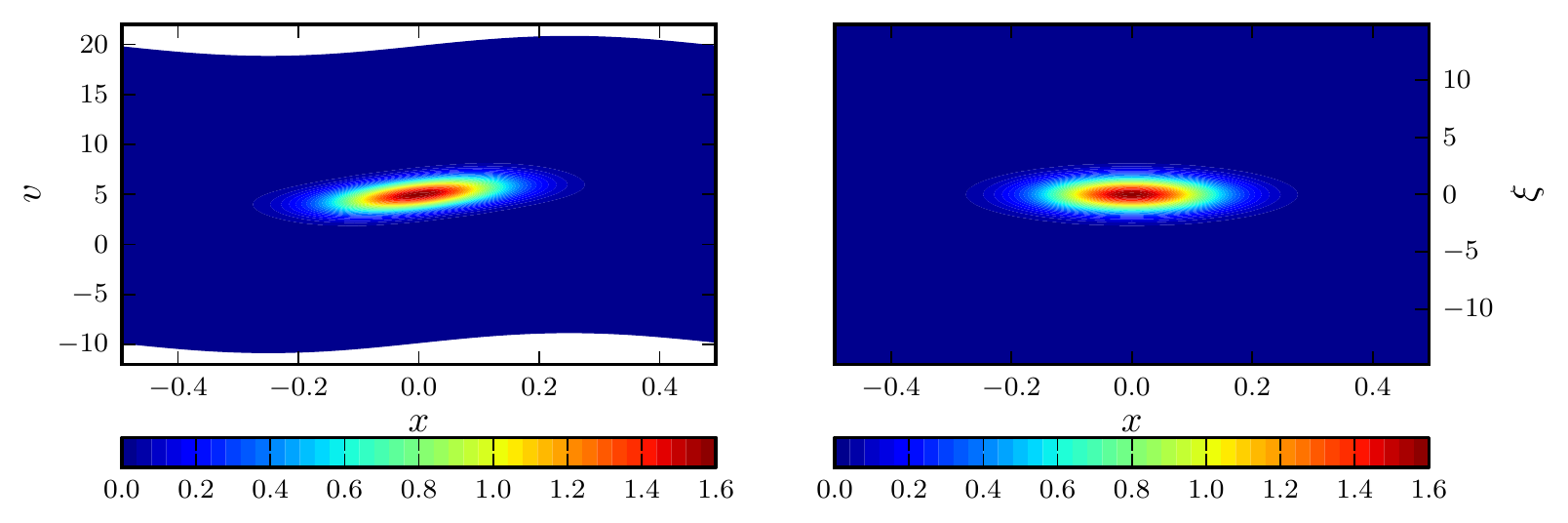}  \\ 
          \includegraphics[scale=1]{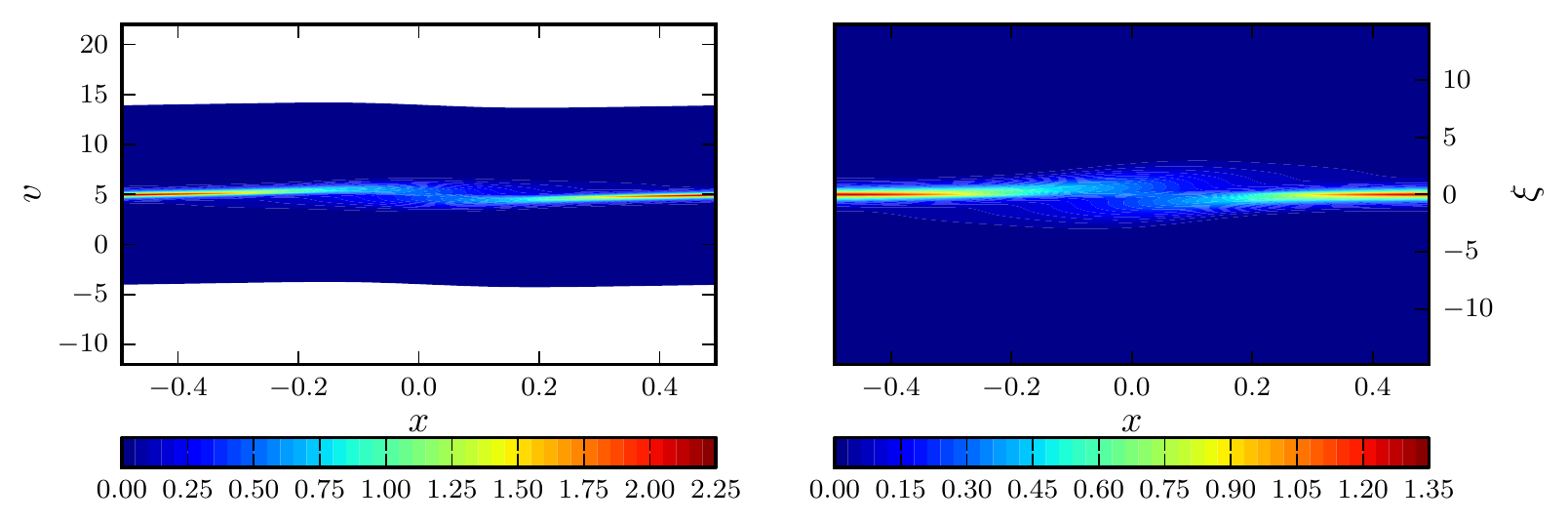}  \\
          \includegraphics[scale=1]{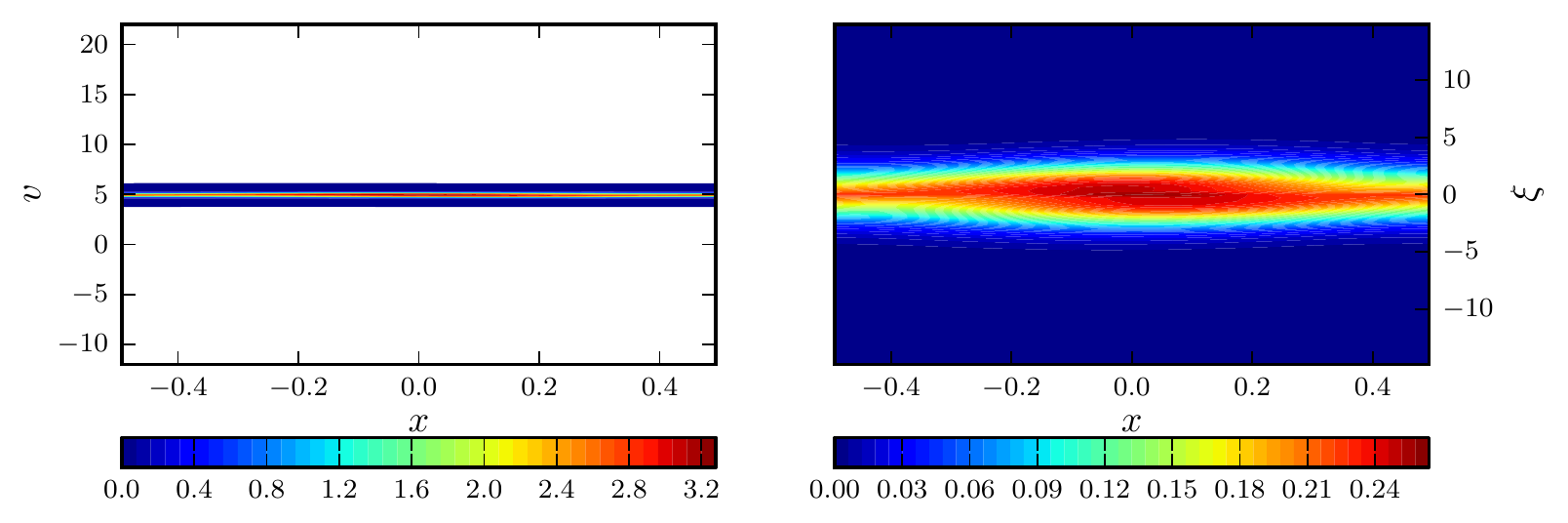} \\
          \includegraphics[scale=1]{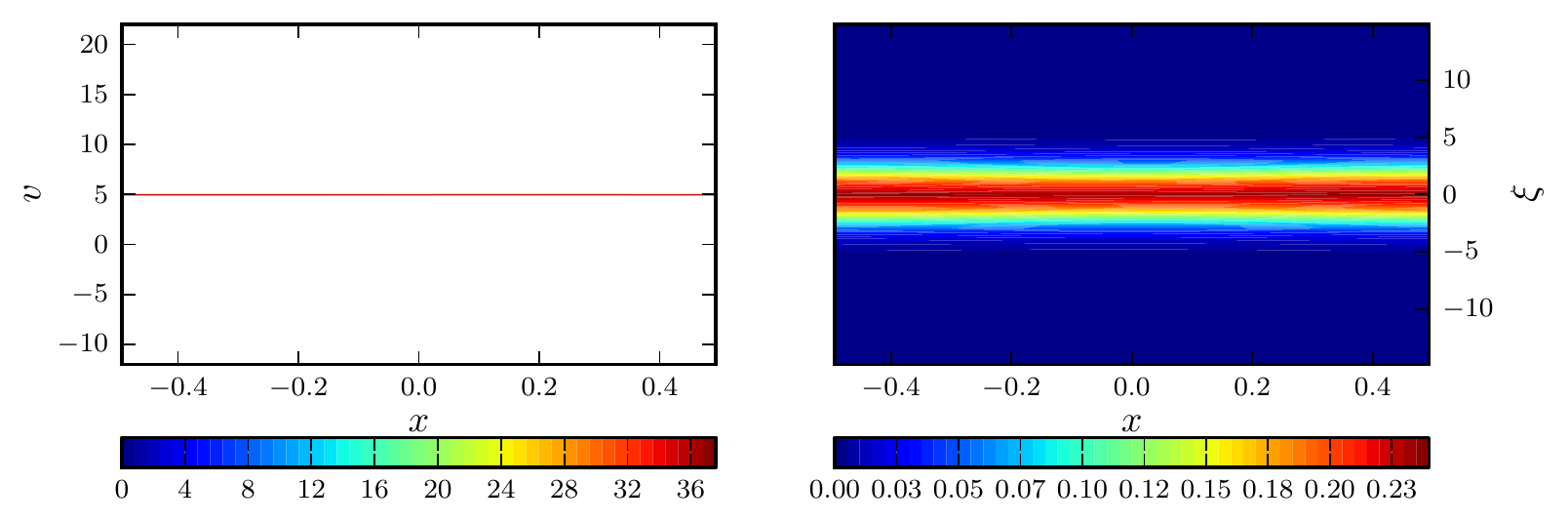} 
        \end{tabular}
        \caption{\textbf{Test 2 -} Contour plot of the original distribution $f(t,x,v)$ (left) and its rescaled counterpart $g(t,x,\xi)$ (right), at times $t=0$, $t=0.5$, $t=2.5$ and $t=5$, in the Motsch-Tadmor case.}
        \label{fig:evolFG_MT}
      \end{figure}
      
     As presented in Figure \ref{fig:evolFG_MT}, on the one hand the distribution $f$ in classical variables $v$ concentrates in
     velocity direction as time evolves. On the other hand, the distribution $g$ after 
     scaling behaves nicely in large time, with neither concentration
     or spreading in $\xi$. As we consider the case where $x$ lies on a torus,
     the equation converges to a global equilibrium, even if the  influence function is local

\begin{figure}
        \includegraphics[scale=1]{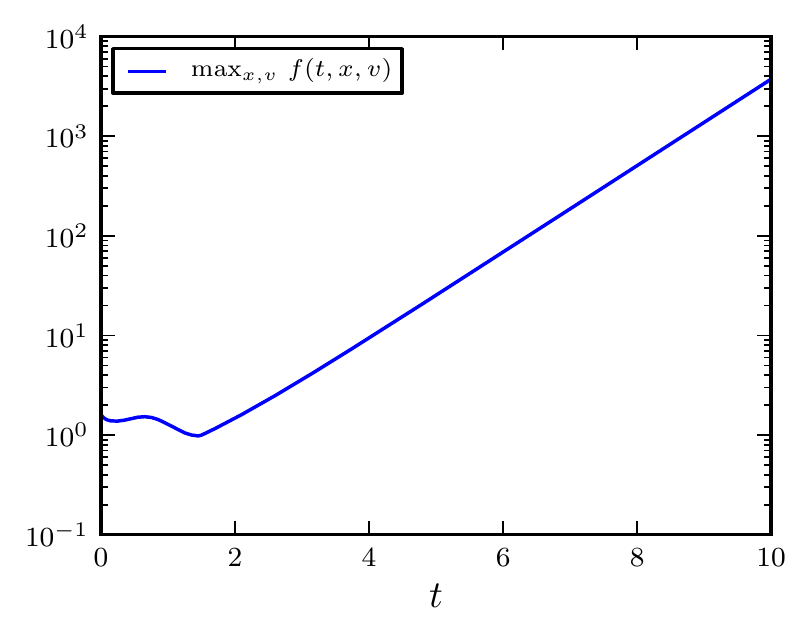}
        \caption{\textbf{Test 2 -} Time evolution of the maximum value of $f$, in the Motsch-Tadmor case.}
        \label{fig:Max_MT}
      \end{figure}

     To further understand the rate of concentration, let us look at
     the maximum value of the reconstructed $f$ against time in Figure \ref{fig:Max_MT}. 
     This quantity will give us a good information on the rate of convergence of $f$ toward a monokinetic distribution.
     We observe an exponential growth in time of this quantity, as expected from the theoretical behavior given by Theorem \ref{thm:flocking}.
     
     \begin{remark}
       Another property of our model is its efficiency, when compared to the original equation. Indeed, the flocking operator for $f$ is written as a convolution in both space and velocity variables, and the numerical cost for its computation with a simple quadrature rule is then proportional to $\mathcal O\left (N_x^2 \, N_\xi^2\right )$.
       The rescaled model, although given by a system of equation, is only obtained thanks to a convolution in space. Its numerical complexity is then proportional to $\mathcal O\left (N_x^2\right )$, which is a huge improvement, specially in higher dimension.
     \end{remark}

    \subsection{Test 3 - One Dimensional Cucker-Smale}
    
      We are finally interested in numerical simulations of the flocking equation \eqref{eq:flock} in the Cucker-Smale case \eqref{eq:collisionMT} for $d=1$ with this time the global influence function
      \[
        \phi(r) = \frac1{(1+r)^{-1/2}},
      \]
      and periodic boundary conditions in $x$.
      We will use for this the rescaled model \eqref{eq:g}. 
      We recall that in the particular Cucker-Smale case, the equation describing the evolution for $g$ is reduced to \eqref{eq:gCS1D}, namely it has one more transport term than the  Motsch-Tadmor model. 
      Moreover, this time, the evolution of $\omega$ is given by the solution to the partial differential equation \eqref{eq:omega}, and is no longer explicit.
      
      We take as an initial condition a step function in the phase space:
      \[
        f_0(x,v) = \bm{1}_{|x| \leq 1/4}(x) \, \bm{1}_{|v| \leq 2}(v), \ \forall x \in \TT, \ v \in \RR.
      \]
      The discretization of equation \eqref{eq:gCS1D} is done as described in section \ref{sub:flockMod}, and we take $N_x = 75$ points in the physical space,  and $N_\xi = 75$ points in the rescaled velocity space for a rescaled velocity variable $\xi \in [-10,10]$. We choose $\Delta t = 1/1200$.
      
      \begin{figure}
        \begin{tabular}{c}
          \includegraphics[scale=1]{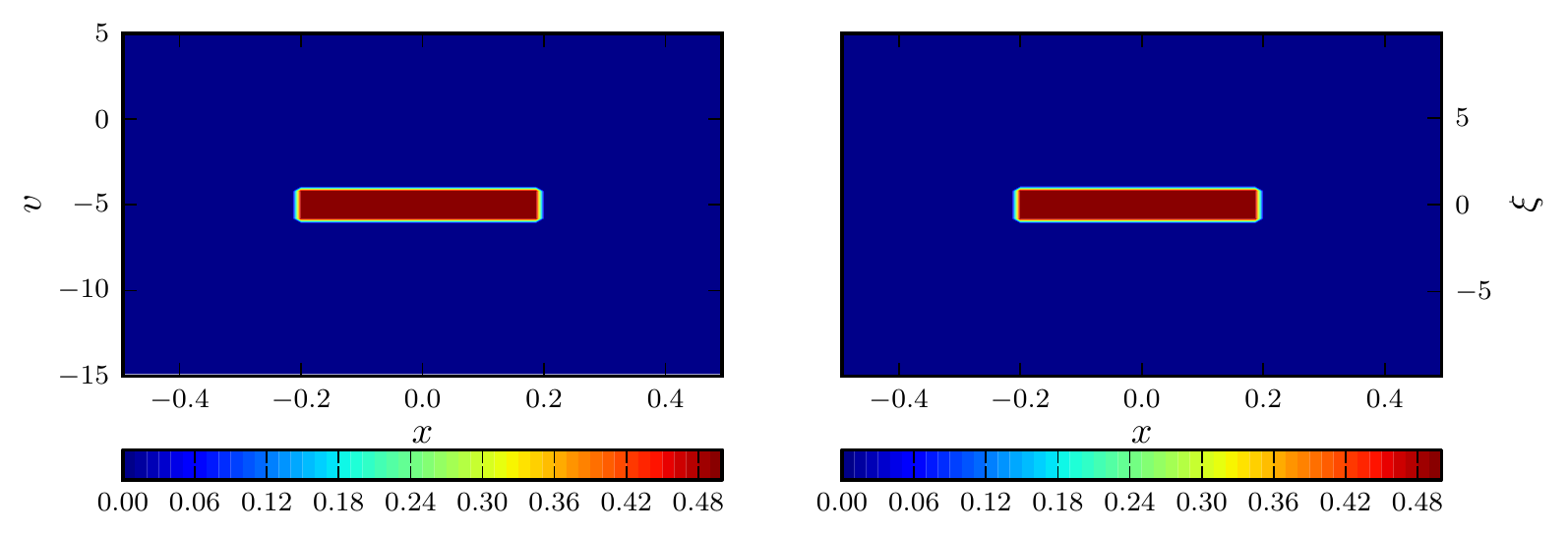}  \\ 
          \includegraphics[scale=1]{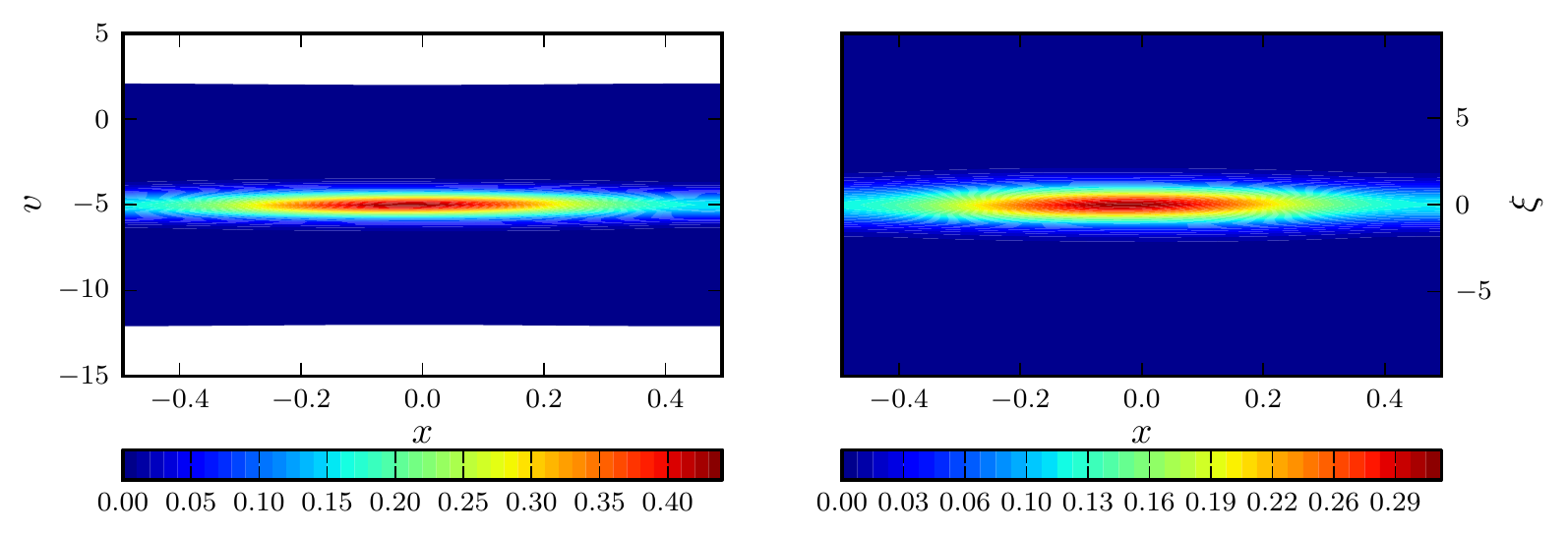}  \\
          \includegraphics[scale=1]{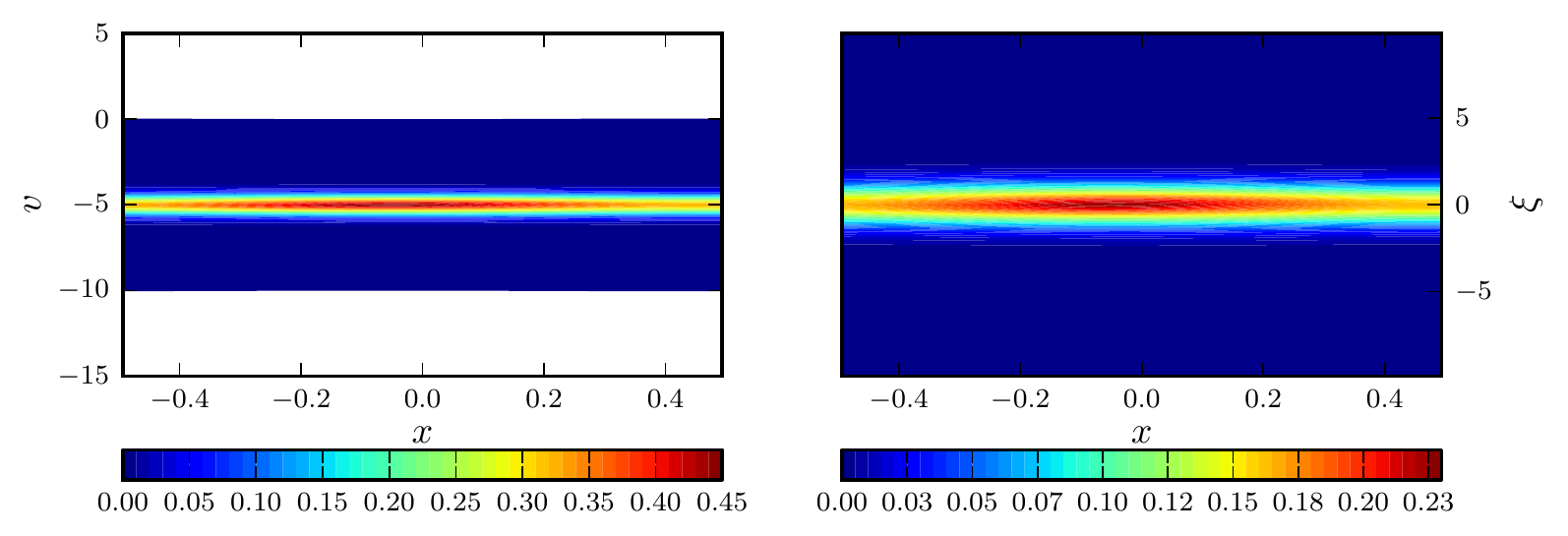} \\
          \includegraphics[scale=1]{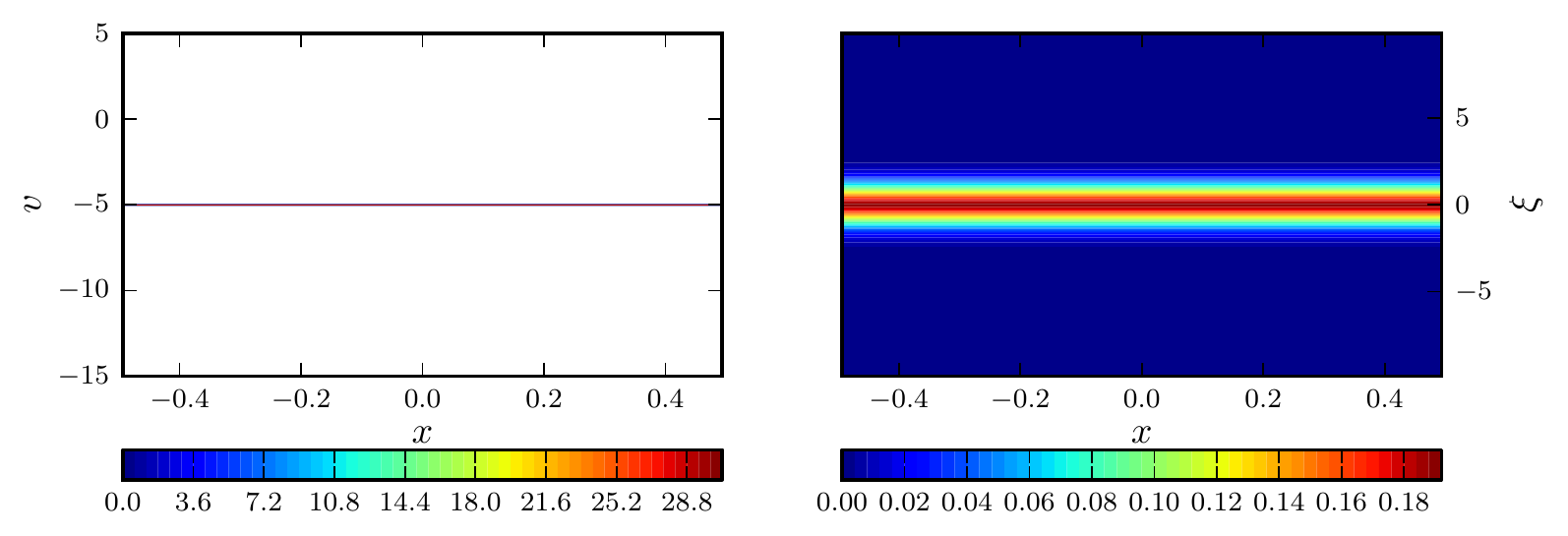} 
        \end{tabular}
        \caption{\textbf{Test 3 -} Contour plot of the original distribution $f(t,x,v)$ (left) and its rescaled counterpart $g(t,x,\xi)$ (right), at times $t=0$, $t=1$, $t=2$ and $t=15$, in the Cucker-Smale case.}
        \label{fig:evolFG_CS}
      \end{figure}

\begin{figure}
        \includegraphics[scale=1]{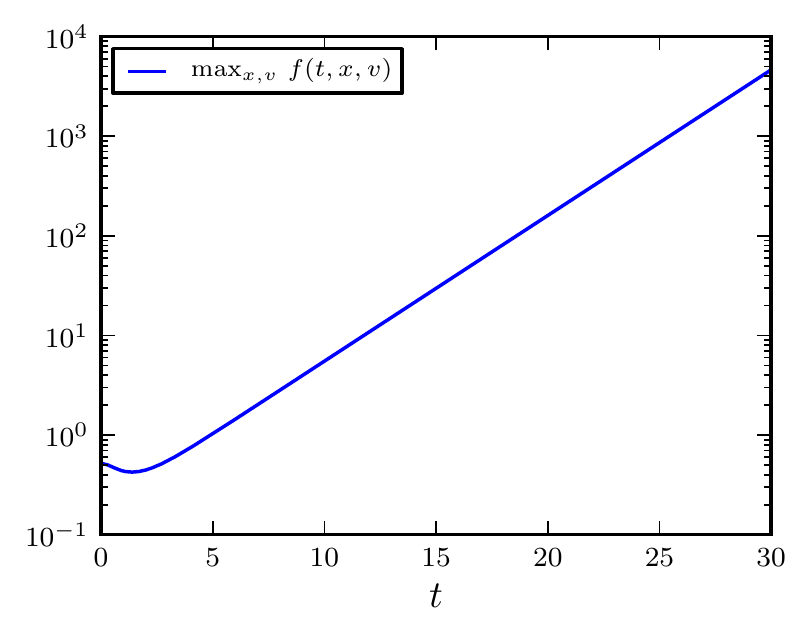}
        \caption{\textbf{Test 3 -} Time evolution of the maximum value of $f$, in the Cucker-Smale case.}
        \label{fig:Max_CS}
      \end{figure}

      We observe in Figure \ref{fig:evolFG_CS} that although this initial condition is not very regular, our first order schemes are dissipative enough to deal with it quite easily.
      More importantly, due to the  Cucker-Smale type of interaction,  particles which are very far from the rest of the flock still have some influence, so the flocking dynamics is slower than the Motsch-Tadmor case, as seen in Figure \ref{fig:Max_CS}.
      Nevertheless, we still have an exponential convergence toward this flock, as predicted by Theorem \ref{thm:flocking}.

  \section*{Acknowledgments}
    Part of this research was conducted during the post-doctoral stay of the first author Thomas Rey (TR) at CSCAMM in the university of Maryland, College Park, under the supervision of Eitan Tadmor. 
    TR would like to warmly thank Eitan Tadmor and all the staff of
    CSCAMM, along with the KI-Net program, for their kindness, their
    welcoming attitude, their availability and by the overall quality
    of his stay.
    The second author Changhui Tan (CT) would like to thank the KI-Net
    program, and in particular, Eitan Tadmor for his consistent help
    and care.
    The research of TR and CT was granted by the NSF Grants DMS
    10-08397, RNMS 11-07444 (KI-Net) and ONR grant N00014-1210318.
    
    Both authors would like to thanks Eitan Tadmor for the very fruitful discussions they had about this work.

\bibliographystyle{acm}
\bibliography{biblioRT}

\end{document}